# THE DISCRETE SCHWARZ-PICK LEMMA FOR CIRCLE PACKINGS REVISITED

ARHAM LODHA

ABSTRACT. The Discrete Schwarz-Pick Lemma is a discrete analogue of the classical result from complex analysis, arising from the connection between circle packings and conformal maps established by Thurston. Previous works by Beardon-Stephanson and Van Eeuwen proved this lemma for circle packings where circles are tangent or intersect at non-obtuse angles, corresponding to inversive distances $I \in [0, 1]$. This paper extends the investigation to circle packings with obtuse intersections ($I \in (-1, 0)$) and disjoint packings ($I > 1$). We prove that the Discrete Schwarz-Pick Lemma holds for the full range of intersecting circle packings with inversive distances in $(-1, 1]$, provided an additional condition on the weights of each triangle is satisfied. The proof relies on a variational principle for circle packings with inversive distances. Conversely, we show that the lemma fails for disjoint circle packings where I≥1. This is demonstrated by constructing a specific counterexample on a triangulated disk with four vertices.

## 1. INTRODUCTION

### 1.1. Background.

In March of 1985, in the International Symposium in Celebration of the Bieberbach Conjecture, William Thurston suggested the connection between circle packings and conformal maps. William Thurston conjectured that the conformal map from a simply connected domain $\Omega$ to the unit disk $\mathbb{D}$ could be approximated by circle packings, a claim famously proven by Rodin and Sullivan in 1987 [1].

Furthermore, in his work constructing hyperbolic structures on 3 manifolds, Thurston introduced the notion of circle packing metrics on triangulated surfaces with non-obtuse intersection angles [2]. The prescription of angles corresponded to the fact that the angles of intersecting circles are invariant under Möbius Transformations. For triangulated surfaces with Thurston's construction there are singularities at the vertices. Classic discrete curvature $K(v)$, which is defined as the angle deficit at vertex $v$, is used to describe the singularity at the vertex.

As a generalization of Thurston's Circle packings, inversive distance circle packings were introduced by Bowers and Stephanson [3]. This generalization allowed for circles of adjacent vertices to be disjoint with the distance between the two circles measured by inversive distance. For more information, see Bowers-Hurdal [4], Stephanson [5], and Guo [6] for more information.

Starting from Thurston's initial idea of approximating conformal maps by discretizing our space with teh use of circle packings and the eventual proof by Rodin and Sullivan, work began on discretizing the ideas and theorems of complex analysis and Riemann Surfaces with the use of circle packings. Important theorems such as the Riemann Mapping Theorem and Uniformization got their discretized equivalents [7].

The theorem of interest for this paper, is the Discrete Schwarz-Pick Lemma, which is the discretized version. of the Schwarz-Pick Lemma in complex analysis.

Beardon and Stephanson established the Discrete Schwarz-Pick Lemma for the case of tangent circle packings using discrete Perron's method [8]. The Discrete Schwarz-Pick Lemma, has been useful in proving other discrete results, for example the uniqueness of an infinite circle packing on the hyperbolic plane. Later, Van Eeuwen proved the case for intersection angles $\varphi$ between two adjacent circles between 0 and $\frac{\pi}{2}$ or inversive distance between 0 and 1 [9]. The difficulty in extending these results to a broader range of inversive distances, lies in the fact that the space of circle packing metrics lacks convexity which breaks many of the arguments needed to establish this result.

This paper relies on a observation by Zhuo [10], to establish some bounds on inversive distances to reestablish convexity to the space of circle packing metrics to prove the Discrete Schwarz-Pick Lemma. The main tools are the variational principle established by Colin de Verdiere [11] and expanded by Guo [6], and foundational results on circle packing geometries by Thurston [2]. Additionally, the paper presents a counter example to the Discrete Schwarz-Pick Lemma for the case of $I \geq 1$ (i.e. tangent or disjoint circles) which was conjectured to hold. This is surprising because of the fact that both infinitesimal [6] and global rigidity [12] were proven for this case. Many of the ideas that go into proving global rigidity in this case are also used to prove the Discrete Schwarz-Pick.

1.2. **Inversive Distance Circle Packings.**

In this subsection, we discuss the concept of inversive distance circle packings introduced by Bowers and Stephanson [3]. For all $n \in \mathbb{N}$, let $[n] := \{1, ..., n\}$.

Let $S$ be a surface with Euclidean or hyperbolic background geometry and a finite simplicial triangulation $\mathcal{T} = \{V, E, F\}$, where $V$, $E$, $F$ represent the vertices, edges, and faces of $\mathcal{T}$. Let $N = |V|$ be a number of vertices. We identify the vertex set $V$ with a index set $[N]$. Thus $E \subset \binom{[N]}{2}$ is the set of edges (pairs of indices) and $F \subset \binom{[N]}{3}$ is the set of faces. We will adopt the following notation convention for convenience:
1. For generic indices represented by variables (e.g. $i, j, k \in [N]$ ) we will denote a edge simply as $ij \in E$ and a face as $ijk \in F$.
2. For specific numerical indices, we will use formal set notation (e.g an edge $\{1, 2\} \in E$ or a face $\{1, 2, 3\} \in F$).

Let $I: E \to (-1, \infty)$ be a function that assigns each edge $\{i, j\}$ a inversive distance $I_{ij} \in (-1, \infty)$. The triple $(S, \mathcal{T}, I)$ will be referred to as a **weighted triangulated surface**.

For all functions $f: V \to \mathbb{R}$, let $f_i := f(i)$ and thus $f$ can be viewed as the vector $(f_1, ..., f_N) \in \mathbb{R}^N$. Let $\mathbb{R}_{>0} := (0, \infty)$. Any map $r: V \to \mathbb{R}_{>0}$ is called a **radius function**. The radius function on a weighted triangulated surface induces a edge length map $l_r: E \to \mathbb{R}_{>0}$ where

$$l_r(ij) = \sqrt{r_i^2 + r_j^2 + 2r_i r_j I_{ij}} \tag{1}$$

if the background geometry is Euclidean and

$$l_r(ij) = \cosh^{-1}\bigl(\cosh(r_i)\cosh(r_j) + I_{ij}\sinh(r_i)\sinh(r_j)\bigr) \tag{2}$$

if the background geometry is hyperbolic. Note that the length $l_r(ij)$ in (1) and (2) is well defined if $r_i, r_j > 0$ and $I_{ij} > -1$. Furthermore the inversive distance $I_{ij}$ determines the geometry of circles attached to $v_i$ and $v_j$ respectively :

- If $I_{ij} \in (-1, 1]$ the circles corresponding to vertices $i$ and $j$ respectively intersect. There exists a unique angle $\Phi_{ij} \in [0, \pi)$ such that $I_{ij} = \cos(\Phi_{ij})$. $\Phi_{ij}$ is the angle formed by the tangent lines of the circles at the intersection point. See Figure 1 for reference. The intersection is obtuse if $I_{ij} \in (-1, 0)$, orthogonal if $I_{ij} = 0$, and acute if $I_{ij} \in (01)$.
- If $I_{ij} = 1$, the circles are tangent.
- If $I_{ij} > 1$, the circles are disjoint.

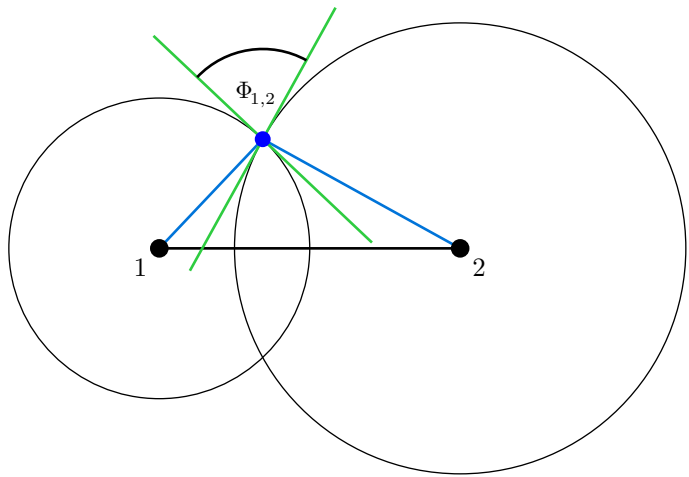

Figure 1. Intersection Angle in Circle Packings

If for each face $\{i, j, k\} \in F$ the induced lengths $l_r$ satisfy the triangle inequalities

$$\begin{aligned} l_r(ij) + l_r(jk) &> l_r(ik) \\ l_r(ik) + l_r(jk) &> l_r(ij) \\ l_r(ij) + l_r(ik) &> l_r(jk) \end{aligned} \quad (3)$$

then $r$ is called a **circle packing metric**. Such a metric induces a combinatorial curvature on the surface $K_r := V \to (-\infty, 2\pi)$ where

$$K_r(i) := 2\pi - \sum_{\alpha > i} \alpha \quad (4)$$

where $\alpha > v$ means all the angles incident to $i$. Note that $r \to K_r$ is a continous function.

### 1.3. **Plan of the Paper.**

The paper is organized as follows. In Section 2 we prove the Discrete Schwarz-Pick Lemma for Circle Packings for a broader range of intersection behaviors for hyperbolic and euclidean background. In Section 3 we present a counterexample to the Discrete Schwarz-Pick Conjecture for the case of tangent and disjoint circles.

## 2. Discrete Schwarz-Pick Lemma for Circle Packing with Intersections

### 2.1. Triangle Preliminaries in Euclidean background geometry.

Let $(S, \mathcal{T}, I)$ be a weighted triangulated surface with Euclidean background geometry and let $r$ be a radius function, and $l_r$ and $K_r$ be the induced edge lengths and discrete curvatures respectively. Suppose $\tau = ijk \in F$, for this section we will temporarily adopt standard geometric conventions where $l_i = l_r(jk)$ and $I_i = I(jk)$ are the edge length and inversive distance of the edge opposite vertex $i$. Let

$$\gamma_\tau = (I_i + I_j I_k, I_j + I_i I_k, I_k + I_i I_j) = (\gamma_\tau^i, \gamma_\tau^j, \gamma_\tau^k). \tag{5}$$

As noted in [6], the edge lengths $l_i, l_j, l_k$ satisfy the triangle inequalities if and only if

$$\begin{aligned}
0 &< (l_i + l_j + l_k)(l_j + l_k - l_i)(l_i + l_k - l_j)(l_i + l_j - l_k) \\
&= \left((l_j + l_k)^2 - l_i^2\right)(l_i - (l_j - l_k))(l_i + (l_j - l_k)) \\
&= \left((l_j + l_k)^2 - l_i^2\right)\left(l_i^2 - (l_j - l_k)^2\right) \\
&= (-l_i^2 + l_j^2 + 2l_j l_k + l_k^2)(l_i^2 - l_j^2 + 2l_j l_k - l_k^2) \\
&= -l_i^4 - l_j^4 - l_k^4 + 2l_i^2 l_j^2 + 2l_i^2 l_k^2 + 2l_j^2 l_k^2
\end{aligned} \tag{6}$$

Substituting the definition for edge lengths defined in (1), by direct computation we have

$$\begin{aligned}
&\frac{1}{4}(l_i + l_j + l_k)(l_j + l_k - l_i)(l_i + l_k - l_j)(l_i + l_j - l_k) \\
&= r_i^2 r_j^2 (1 - I_k^2) + r_i r_k (1 - I_j^2) + r_j^2 r_k^2 (1 - I_i^2) \\
&\quad + 2 r_i r_j r_k (r_i \gamma_\tau^i + r_j \gamma_\tau^j + r_k \gamma_\tau^k)
\end{aligned} \tag{7}$$

Thus we have the following result on the Euclidean triangle inequalities.

**Lemma 2.1.1** *([13], Lemma 2.1)* Suppose $(S, \mathcal{T}, I)$ is a weighted triangulated surface, $r$ is a radius function, and $\tau = \triangle v_i v_j v_k$ is a topological triangle in $F$. Then the edge lengths defined in (1) satisfy the triangle inequalities if and only if

$$\begin{aligned}
r_i^2 r_j^2 (1 - I_k^2) + r_i r_k (1 - I_j^2) + r_j^2 r_k^2 (1 - I_i^2) \\
+ 2 r_i r_j r_k (r_i \gamma_\tau^i + r_j \gamma_\tau^j + r_k \gamma_\tau^k) > 0
\end{aligned} \tag{8}$$

The following corollary obtained in [10] directly follows from `Lemma 2.1.1`.

**Corollary 2.1.2** If $I_i, I_j, I_k \in (-1, 1]$ and $\gamma_\tau \in [0, \infty)^3$, then the triangle inequalities are satisfy for any $(r_i, r_j, r_k) \in \mathbb{R}^3_{>0}$.

It immediately follows that if for all $\tau \in F$ we have that $\gamma^\tau \in \mathbb{R}^3_{\geq 0}$, then the space of circle packing metrics on $(S, \mathcal{T}, I)$ is $\mathbb{R}^N_{>0}$.

### 2.2. Infinitesimal Rigidity of Euclidean Inversive Distance Circle Packings.

Let $(S, \mathcal{T}, I)$ be a weighted triangulated surface with euclidean background where $-1 < I \leq 1$. Let $\tau = ijk \in F$ where $\gamma^\tau \in \mathbb{R}^3_{\geq 0}$. By `Corollary 2.1.2`, $\forall (r_i, r_j, r_k) \in \mathbb{R}^3_{>0}$ induces edge lengths which form a triangle. Let $(u_i, u_j, u_k) = (\ln(r_i), \ln(r_j), \ln(r_k))$. For all $a \in \{i, j, k\}$, let $\theta_a$ inner angle of $\tau$ at the vertex $a$. We have the following useful lemma regarding the monotonicity of $\theta_a$ with respect to $(u_i, u_j, u_k)$, and thus $(r_i, r_j, r_k)$.

**Lemma 2.2.3** *([13])* For all $\tau = ijk \in F$ where $\gamma_\tau \in \mathbb{R}^3_{\geq 0}$, we have the following:

1. $\frac{\partial \theta_a}{\partial u_b} = \frac{\partial \theta_b}{\partial u_a} > 0$ for $a, b \in \{i, j, k\}$ where $a \neq b$.
2. $\frac{\partial \theta_a}{\partial u_a} < 0$ for $a \in \{i, j, k\}$ .
3. The jacobian $J_u(\theta)$ is symmetric and negative semi definite, with null space

$$\mathbb{R}[1, 1, 1] = \{[t, t, t] \mid t \in \mathbb{R}\}. \tag{9}$$

4. The differential 1-form $\sum_{a \in \{i,j,k\}} \theta_a \, du_a$ is closed on $\mathbb{R}^3$ (and hence exact) and

$$W_\tau(u) = \int_0^u \sum_{a \in \{i,j,k\}} \theta_a \, du_a \tag{10}$$

is a well defined concave function of $u \in \mathbb{R}^3$ satisfying

$$\nabla W_\tau = (\theta_i, \theta_j, \theta_k) \tag{11}$$

Furthermore, $W_\tau$ is strictly concave when restricted to $\{tu + (1-t)\hat{u} \mid t \in \mathbb{R}\}$, provided that $u - \hat{u} \neq (a, a, a)$ for $a \in \mathbb{R}$.

**Proof** Part 1 of the lemma comes from Lemma 2.5 in Xu [13]. Part 2 holds because $\theta_i + \theta_j + \theta_k = \pi$. Thus

$$\frac{\partial \theta_i}{\partial u_i} + \frac{\partial \theta_j}{\partial u_i} + \frac{\partial \theta_k}{\partial u_i} = 0. \tag{12}$$

Since $\frac{\partial \theta_a}{\partial u_b} > 0$ for $a \neq b$, by part 1 we deduce that $\frac{\partial \theta_i}{\partial u_i} < 0$. Part 3 and 4 hold because of Lemma 2.6 and Lemma 2.10 respectively from Xu [13]. □

2.3. **Triangle Preliminaries in hyperbolic background geometry.**

Let $(S, \mathcal{T}, I)$ be a weighted triangulated surface with hyperbolic background geometry and let $r$ be a radius function, and $l_r$ and $K_r$ be the induced edge lengths and discrete curvatures respectively. Suppose $\tau = ijk \in F$, for this section we will temporarily adopt standard geometric conventions where $l_i = l_r(jk)$ and $I_i = I(jk)$ are the edge length and inversive distance of the edge opposite vertex $i$. Let

$$\gamma_\tau = (I_i + I_j I_k, I_j + I_i I_k, I_k + I_i I_j) = (\gamma_\tau^i, \gamma_\tau^j, \gamma_\tau^k). \tag{13}$$

Let

$$S_i = \sinh(r_i) \text{ and } C_i = \cosh(r_i). \tag{14}$$

Working the algebra to a inequality similar to (6), we get the following lemma regarding the circle packing metrics for our triangle $\tau$.

**Lemma 2.3.4** *([6], Lemma 3.1)* Suppose $(S, \mathcal{T}, I)$ is a weighted triangulated surface with hyperbolic background geometry. Let $r$ be be a radius function and $\tau = ijk \in F$. Then the edge lengths defined in (1) satisfy the triangle inequalities if and only if

$$2S_i^2 S_j^2 S_k^2 (1 + I_i I_j I_k) + S_i^2 S_j^2 (1 - I_k^2) + S_i^2 S_k^2 (1 - I_j^2) + S_j^2 S_k^2 (1 - I_i^2)$$
$$+ 2S_i S_j S_k (S_i C_j C_k \gamma_\tau^i + C_i S_j C_k \gamma_\tau^j + C_i C_j S_k \gamma_\tau^k) > 0 \tag{15}$$

The following corollary obtained in [10] directly follows from `Lemma 2.3.4`.

**Corollary 2.3.5** *([13])*   If $I_i, I_j, I_k \in (-1, 1]$ and $\gamma_\tau \in [0, \infty)^3$, then the triangle inequalities are satisfy for any $(r_i, r_j, r_k) \in \mathbb{R}^3_{>0}$. Thus $\Omega_\tau = \mathbb{R}^3_{>0}$.

2.4. **Infinitesimal Rigidity of Hyperbolic Inversive Distance Circle Packings.**

Let $(S, \mathcal{T}, I)$ be a weighted triangulated surface with hyperbolic background geometry where $-1 < I \leq 1$. Let $\tau = ijk \in F$ where $\gamma^\tau \in \mathbb{R}^3_{\geq 0}$. By `Corollary 2.3.5`, $\forall (r_i, r_j, r_k) \in \mathbb{R}^3_{>0}$ induces edge lengths which form a triangle. Let

$$u = (u_i, u_j, u_k) = \left( \ln\left(\tanh\left(\frac{r_i}{2}\right)\right), \ln\left(\tanh\left(\frac{r_j}{2}\right)\right), \ln\left(\tanh\left(\frac{r_k}{2}\right)\right) \right) \tag{16}$$

. Let $\theta_a$ inner angle of $\tau$ at the vertex $v_a$ for $a \in \{i, j, k\}$. We have the following useful lemma regarding the monotonicity of $\theta_a$ with respect to $(u_i, u_j, u_k)$, and thus $(r_i, r_j, r_k)$.

**Lemma 2.4.6** *([13])*   For all $\tau = ijk \in F$ where $\gamma_\tau \in \mathbb{R}^3_{\geq 0}$, we have the following:

1. $\frac{\partial \theta_a}{\partial u_b} = \frac{\partial \theta_b}{\partial u_a} > 0$ for $a, b \in \{i, j, k\}$ where $a \neq b$.
2. $\frac{\partial \theta_a}{\partial u_a} < 0$ for $a \in \{i, j, k\}$ .
3. The jacobian $J_u(\theta)$ is symmetric and negative semi definite, with null space

$$\mathbb{R}[1, 1, 1] = \{[t, t, t] \mid t \in \mathbb{R}\}. \tag{17}$$

4. The differential 1-form $\sum_{a \in \{i,j,k\}} \theta_a \, du_a$ is closed on $\mathbb{R}^3$ (and hence exact) and

$$W_\tau(u) = \int_0^u \sum_{a \in \{i,j,k\}} \theta_a \, du_a \tag{18}$$

is a well defined concave function of $u \in \mathbb{R}^3$ satisfying

$$\nabla W_\tau = (\theta_i, \theta_j, \theta_k) \tag{19}$$

Furthermore, $W_\tau$ is strictly concave when restricted to $\{tu + (1-t)\hat{u} \mid t \in \mathbb{R}\}$, provided that $u - \hat{u} \neq (a, a, a)$ for $a \in \mathbb{R}$.

**Proof**   Part 1 of the lemma holds by Lemma 3.6 in Xu [13]. By Remark 11 in Xu [13],

$$\frac{\partial \theta_i}{\partial u_i} + \frac{\partial \theta_j}{\partial u_i} + \frac{\partial \theta_k}{\partial u_i} < 0 \tag{20}$$

since $\frac{\partial \theta_a}{\partial u_b} > 0$ for $a \neq b$ by 1, $\frac{\partial \theta_a}{\partial u_a} < 0$, thus Part 2 holds. Part 3 holds by Lemma 3.7 in Xu [13]. Finally, Part 3 holds by Lemma 3.9 in Xu [13]. □

2.5. **Degeneration of Circle Packing Metrics.**

Let $(S, \mathcal{T}, I)$ be a closed connected weighted triangulated surface with Euclidean or hyperbolic ground geometry where $-1 < I \le 1$. Let $\Phi : E \to [0, \pi)$ be the corresponding intersection angle for $I$, ie $I(e) = \cos(\Phi(e))$. $\forall J \subset V$, let $F_J$ be the subcomplex of $\mathcal{T}$ consisting of simplicies whose vertices are in $J$ and

$$\operatorname{Lk}(J) = \{(\tau, v) \in F \times J \mid \{v\} = \tau \cap J\}. \tag{21}$$

.

**Lemma 2.5.7** Let $J$ be some subset of $V$. Let $(r_n)_{n \in \mathbb{N}}$ be a sequence of circle packing metrics for $(S, \mathcal{T}, I)$ $(-1 < I \le 1)$ where if $j \in J$

$$\lim_{n \to \infty} r_n(j) = 0 \tag{22}$$

and if $k \notin J$,

$$\lim_{n \to \infty} r_n(k) > 0. \tag{23}$$

Then

$$\lim_{n \to \infty} \sum_{j \in J} K_{r_n}(j) = 2\pi \chi(F_J) - \sum_{(\tau, v) \in \operatorname{Lk}(J)} (\pi - \Phi(\tau \setminus \{v\})) \tag{24}$$

**Proof** This proof follows the argument outlined in Proposition 4.1 in Chow and Luo [14]. $\forall \tau \in F$ and $\forall i \in \tau$, $\theta_i^\tau$ be the inner angle of $i$ in triangle $\tau$. Let $T_1$, $T_2$, and $T_3$ be the sets of triangles with one, two, and three vertices in $J$. Thus for any circle packing metric $r$, we have following:

$$\begin{aligned}
\sum_{j \in J} K_r(j) &= 2\pi |J| - \sum_{j \in J} \sum_{\substack{\tau \in F \\ j \in \tau}} \theta_j^\tau \\
&= 2\pi |J| - \sum_{\tau \in F} \sum_{j \in \tau \cap J} \theta_j^\tau \\
&= 2\pi |J| - \left( \sum_{\substack{\tau \in T_1 \\ i \in \tau \cap J}} \theta_i^\tau + \sum_{\substack{\tau \in T_2 \\ i, j \in \tau \cap J \text{ distinct}}} \theta_i^\tau + \theta_j^\tau + \sum_{\substack{\tau \in T_3 \\ \tau = ijk}} \theta_i^\tau + \theta_j^\tau + \theta_k^\tau \right)
\end{aligned} \tag{25}$$

Note there is bijection between $T_1$ and $\operatorname{Lk}(J)$, so we have the following equivalence:

$$\sum_{\substack{\tau \in T_1 \\ i \in \tau \cap J}} \theta_i^\tau = \sum_{(i, \tau) \in \operatorname{Lk}(J)} \theta_i^\tau \tag{26}$$

Given that $\forall i \in J$, $r_n(i) \to 0$ we know the following:

1. $\forall \tau = ijk \in T_3$ as $(r_n(i), r_n(j), r_n(k)) \to (0, 0, 0)$, $ijk$ approaches a Euclidean triangle. Thus $\theta_i^\tau + \theta_j^\tau + \theta_k^\tau \to \pi$.
2. $\forall \tau \in T_2$ and $i, j \in \tau \cap J$ distinct as $(r_n(i), r_n(j)) \to (0, 0)$, $ijk$ approaches a geodesic segment. Thus $\theta_k^\tau \to 0 \Rightarrow \theta_i^\tau + \theta_j^\tau \to \pi$.
3. $\forall \tau = ijk \in T_1$ and $i \in \tau \cap J$, then $i$ approaches the intersection point of the circles corresponding to $j$ and $k$. Thus $\theta_i^\tau \to \pi - \Phi(jk) = \pi - \Phi(\tau \setminus \{i\})$.

Putting all that together, we see:

$$\lim_{n\to\infty}\sum_{j\in J} K_{r_n}(j) = 2\pi|J| - \left(\sum_{(i,\tau)\in \mathrm{Lk}(J)} \pi - \Phi(\tau\setminus\{i\}) + |T_2|\pi + |T_3|\pi\right)$$

$$= 2\pi\left(|J| - \frac{1}{2}|T_2| - \frac{1}{2}|T_3|\right) - \sum_{(i,\tau)\in \mathrm{Lk}(J)} \pi - \Phi(\tau\setminus\{i\}) \quad (27)$$

$$= 2\pi\left(|J| - \frac{1}{2}(|T_2| + 3|T_3|) + |T_3|\right) - \sum_{(i,\tau)\in \mathrm{Lk}(J)} \pi - \Phi(\tau\setminus\{i\})$$

$|J|$ is the number of vertices in $F_J$. $|T_3|$ is the number of 2 simplicies in $F_J$. By construction, the number of edges is $\frac{1}{2}(|T_2| + 3|T_3|)$. Thus

$$\lim_{n\to\infty}\sum_{j\in J} K_{r_n}(j) = 2\pi\chi(F_J) - \sum_{(\tau,v)\in \mathrm{Lk}(J)} (\pi - \Phi(\tau\setminus\{v\})) \quad (28)$$

□

We get the following corollary from the lemma above.

**Corollary 2.5.8** If $\forall \tau \in F$ we have that $\gamma^\tau \in \mathbb{R}^3_{\geq 0}$, then $\forall r: V \to \mathbb{R}^3_{>0}$ circle packing metrics

$$\sum_{j\in J} K_r(v_j) > 2\pi\chi(F_J) - \sum_{(\tau,v)\in \mathrm{Lk}(J)} (\pi - \Phi(\tau\setminus\{v\})) \quad (29)$$

**Proof** Since the background geometry is Euclidean or Hyperbolic, we know that the maximum sum of interior angles of a triangle is $\pi$. Thus for any circle packing metric $r$,

$$\sum_{\substack{\tau\in A_2 \\ i,j\in\tau\cap J \text{ distinct}}} \theta^\tau_i + \theta^\tau_j < \pi|T_2| \text{ and } \sum_{\substack{\tau\in A_3 \\ i,j,k\in\tau\cap J \text{ distinct}}} \theta^\tau_i + \theta^\tau_j + \theta^\tau_k < \pi|T_3|. \quad (30)$$

Furthermore since $\forall \tau \in F$ we have that $\gamma^\tau \in \mathbb{R}_{>0}$, we can apply `Lemma 2.2.3` and `Lemma 2.4.6` to get monotonicity of interior angles. $\forall \tau \in A_1$ where $i \in \tau \cap J$, we know that $\theta^\tau_i$ is monotonically increasing, thus $\pi - \Phi(\tau\setminus\{v\})$ is the upper bound for $\theta^\tau_i$. Thus we have:

$$\sum_{j\in J} K_r(j) > 2\pi|J| - \left(\sum_{(i,\tau)\in \mathrm{Lk}(J)} \pi - \Phi(\tau\setminus\{i\}) + |T_2|\pi + |T_3|\pi\right) \quad (31)$$

Proceeding similarly to `Lemma 2.5.7`, we see that

$$\sum_{j\in J} K_r(v_j) > 2\pi\chi(F_J) - \sum_{(\tau,v)\in \mathrm{Lk}(J)} (\pi - \Phi(\tau\setminus\{v\})) \quad (32)$$

□

2.6. **Main Theorem.**

**Theorem 2.6.9** Let $(S, \mathcal{T}, I)$ be a closed connected weighted triangulated surface where $I \in (-1, 1]$. Let $\Phi$ be its corresponding intersection angle function. Let $A, B \subset [N]$ where $N = A \sqcup B$ and $A \neq \varnothing$ and in the case Euclidean background, $B \neq \varnothing$. If $\forall \tau \in F$ we have that $\gamma^\tau \in \mathbb{R}^3_{\geq 0}$, then for any circle packing metrics $R$ and $r$ for $(S, \mathcal{T}, I)$ where $R|_B \geq r|_B$ and $K_R|_A \geq K_r|_A$ it means that $R \geq r$.

**Proof**

Let
$$X := \{x \in \mathbb{R}^N_{>0} \mid x|_B \geq x|_B \text{ and } K_x|_A \geq K_r|_A\}. \tag{33}$$

Note that $r, R \in X$ by assumption. We want to show that $\forall v \in V$,
$$r(v) = \inf\{x(v) : x \in X\}. \tag{34}$$

This will imply that $R \geq r$ and finish the proof.

$\forall a, b \in X$, let $c_i = \min(a_i, b_i)$, thus $c = \min(a, b) := (c_1, ..., c_N)$. Clearly, $c_j \geq r_j$ for all $j \in B$. $\forall v \in V_A$, let $\{u_1, ..., u_m\}$ the set of vertices adjacent to $v$. Assume without a loss of generality that, $c(v) = a(v)$. By assumption $c(u_i) \leq a(u_i)$, thus by part 1 of `Lemma 2.2.3` and part 1 of `Lemma 2.4.6`, $K_c(v) \geq K_a(v)$. Thus $K_c(v) \geq K_a(v) \geq K_r(v)$. Thus $c \in X$ and $X$ is closed under minimums.

Let $w = \inf X$ (component-wise). $\forall v \in V$, let $(x^v_n)_{n \in \mathbb{N}}$ be a sequence in $X$ such that
$$w(v) = \lim_{n \to \infty} x^v_n(v). \tag{35}$$

Let $(y_n)_{n \in \mathbb{N}}$ be a sequence where
$$y_n := \min\{x^v_n \mid v \in V\}. \tag{36}$$

By closure under minimums, $y_n \in X$ for all $n \in \mathbb{N}$. Since $\forall n \in \mathbb{N}$ and $v \in V$ $y_n \leq x^v_n$, that means
$$\limsup_{n \to \infty} y_n(v) \leq \limsup_{n \to \infty} x^v_n(v) = \lim_{n \to \infty} x^v_n(v) = w(v). \tag{37}$$

At the same time since $w$ is the componentwise infimum that means that $\forall n \in \mathbb{N}$ $y_n(v) \geq w(v)$ which means
$$w(v) \leq \liminf_{n \to \infty} y_n(v) \Rightarrow \lim_{n \to \infty} y_n(v) = w(v). \tag{38}$$

We will use this sequence $y_n$ to verify that $w \in X$. We have to check that $w \in \mathbb{R}^N_{>0}$, $w|_B \geq r|_B$, and $K_w|_A \geq K_r|_A$.

Assume the contrary that $w \notin \mathbb{R}^N_{>0}$ thus $\exists v \in V$ such that $w(v) = 0$. Let $J := \{j \in [N] \mid w(v_j) = 0\}$. By `Lemma 2.5.7`, we have that
$$\lim_{n \to \infty} \sum_{i \in J} K_{y_n}(v_i) = 2\pi \chi(F_J) - \sum_{(\tau, v) \in \mathrm{Lk}(J)} (\pi - \Phi(\tau \setminus \{v\})). \tag{39}$$

Note that $\forall v \in B$, $w(v) \geq r(v) > 0 \Rightarrow J \cap B = \varnothing$ which means $J \subset A$. Thus $\forall j \in J$,
$$K_{y_n}(v_j) \geq K_r(v_j) \tag{40}$$

which means
$$\lim_{n \to \infty} \sum_{j \in J} K_{y_n}(v_j) \geq \sum_{j \in J} K_r(v_j). \tag{41}$$

Applying `Corollary 2.5.8` for $\sum_{j \in J} K_r(v_j)$ we see that

$$\lim_{n\to\infty}\sum_{j\in J}K_{y_n}(v_j)\geq\sum_{j\in J}K_r(v_j)>2\pi\chi(F_J)-\sum_{(\tau,v)\in\mathrm{Lk}(J)}(\pi-\Phi(\tau\setminus\{v\}))\quad(42)$$

which is a contradiction. Thus by contradiction $w\in\mathbb{R}_{>0}^N$ and $w$ is a circle packing metric for $(S,\mathcal{T},I)$.

$\forall n\in\mathbb{N}$, we know that $y_n|_B\geq r|_B$ and

$$w=\lim_{n\to\infty}y_n\quad(43)$$

which implies that $w|_B\geq r|_B$. By the definition of infimum, $r|_B\geq w|_B$. Thus $r|_B=w|_B$.

Since $r\to K_r$ is a continuous function, $\forall i\in A$,

$$K_w(v_i)=\lim_{n\to\infty}K_{y_n}(v_i)\quad(44)$$

by the sequential definition of continuity. $K_{y_n}(v_i)\geq K_r(v_i)\Rightarrow K_w(v_i)\geq K_r(v_i)$. Thus $w=X$.

We know that $w=X$ and $w|_B=r|_B$, we just need to show that $w|_A=r|_A$. Assume the contrary that $\exists i\in A$ such that $K_w(v_i)>K_r(v_i)$. $\forall t\in(0,w(v_i))$, let $s_t:V\to\mathbb{R}_{>0}^+$ where

$$s_t(v)=\begin{cases}w(v)-t\text{ if }v=v_i\\w(v)\text{ otherwise}\end{cases}\quad(45)$$

Since $i\in A$, $s_t|_B=w|_B\geq x|_B$. Let $U\subset V$ be the subset of adjacent vertices to $i$. Since $s_t|_U=w|_U$ and $s_t(i)<w(i)$. By Part 1 of both `Lemma 2.2.3` and `Lemma 2.4.6`, we $K_{s_t}(i)<K_w(i)$ and $K_{s_t}(j)>K_w(j)\geq K_r(j)$ for all $j\in U$. By continuity, $\exists t_0\in(0,w(i))$ such that $K_r(i)\leq K_{s_t}(i)<K_w(i)$. Thus $s_{t_0}\in X$, however note that $s_{t_0}(i)<w(i)$. This is a contradiction since $w$ is a componentwise infimum. Thus by contradiction $K_w|_A=K_r|_A$.

Now we have to show that $r=w$ using the facts that $r|_B=w|_B$ and $K_r|_A=K_w|_A$. Let

$$S=\begin{cases}\mathbb{R}\text{ if background geometry is Euclidean}\\(0,\infty)\text{ if background geometry is hyperbolic}\end{cases}\quad(46)$$

Let $\varphi:(0,\infty)\to S$ where

$$\varphi(x):=\begin{cases}\ln(x)\text{ if background geometry is Euclidean}\\\ln(\tanh(\frac{x}{2}))\text{ if background geometry is hyperbolic}\end{cases}\quad(47)$$

$\varphi$ is a homeomorphism so $\varphi^{-1}$ exists and is a continuous. Let $\mathcal{R}:=\{x\in\mathbb{R}_{>0}^N:\varphi^{-1}|_B=r|_B\}$. Let $W:S^N\to\mathbb{R}$ where

$$W(x):=\sum_{\tau=ijk\in F}W_\tau(x_i,x_j,x_k)\quad(48)$$

where $W_\tau$ is defined in `Lemma 2.2.3` and `Lemma 2.4.6`. $\forall i\in V$,

$$\begin{aligned}
\frac{\partial W}{\partial x_a} &= \sum_{\tau=ijk\in F} \frac{\partial W_\tau}{\partial x_a} \\
&= \sum_{\tau=ajk\in F} \frac{\partial W_\tau}{\partial x_a} \\
&= \sum_{\tau=ajk\in F} \theta_a^\tau \\
&= 2\pi - K_{\varphi^{-1}(x)}(a),
\end{aligned} \tag{49}$$

where the second to last equality comes from Part 2 of `Lemma 2.2.3` and `Lemma 2.4.6` and the last equality comes from the definition of discrete curvature. Thus

$$\nabla W(x) = 2\pi - K_{\varphi^{-1}(x)}. \tag{50}$$

Let $F = W|_\mathcal{R}$, thus $F : \mathbb{R}_{>0}^A \to \mathbb{R}$ and

$$\nabla F(x) = 2\pi - K_{\varphi^{-1}(x)}|_A. \tag{51}$$

$\forall x, y \in \mathcal{R}$ such that $xx \neq y$, take $M = \{i \in V \mid x(i) \neq y(i)\}$ and $N = \{i \in V \mid x(i) = y(i)\}$. Note $M \neq \varnothing$, since $x \neq y$, $B \subset N$ and $M \cap N = \varnothing$. Since $(S, \mathcal{T})$ is connected, then the graph of the triangulation is connected. Thus $\exists (i,j) \in M \times N$ such that $\exists k \in V$ such that $ijk \in F$. By definition

$$(x(i) - y(i), x(j) - y(j), x(k) - y(k)) \neq (0,0,0) \tag{52}$$

because $x(i) \neq y(i)$. Thus $W_{ijk}$ is strictly concave on $\{tu + (1+t)\hat{u}\}$ where $u = (x(i), x(j), y(k))$ and $\hat{u} = (y(i), y(j), y(k))$. Thus $W_\tau$ is strictly concave on $\mathcal{R}$, thus $W$ and $F$ are strictly concave on $\mathcal{R}$. Which implies that $\nabla F$ is injective and $x \to K_{\varphi^{-1}(x)}|_A$ is injective. Thus since $K_w|_A = K_r|_A \Rightarrow w|_A = r|_A \Rightarrow w = r$. Thus $r = \inf(X)$ and $R \geq r$. $\square$

## 3. The Discrete Schwarz-Pick Lemma for $I \geq 1$

For $I \geq 1$, it was conjectured that the Discrete Schwarz-Pick Lemma would hold for all $I \geq 1$, because for such circle packings it is known that rigidity holds, for hyperbolic (respectively Euclidean) geometry discrete curvature determines (respectively up to scale) the circle packing metric, by work done by Luo [12]. In such circle packings, circles corresponding to adjacent vertices are allowed to be disjoint or tangent.

**Conjecture 3.10** Let $(S, \mathcal{T}, I)$ be a closed connected weighted triangulated surface in Euclidean or hyperbolic geometry such that $I \geq 1$. Let $V = A \sqcup B$ such that $A \neq \varnothing$ and $B \neq \varnothing$, if the background geometry is Euclidean. Then for all circle packing metrics $R$ and $r$ on $(S, \mathcal{T}, I)$ where $R|_B \geq r|_B$ and $K_R|_A \geq K_r|_A$, $R \geq r$.

### 3.1. **Counterexample for the Conjecture.**

The statement of the conjecture concerns closed surface, but a counterexample on a triangulated disk is sufficient. By a standard doubling argument, any such counterexample can be extended to a closed surface.

Let $(S, \mathcal{T}, I)$ be a disk with boundary with interior vertices $A$ and boundary vertices $B$. Construct a closed surface $S'$ by taking two copies of $S$ and identifying their boundaries with the identity map. This induces a triangulation $\mathcal{T}'$ and inversive distance $I'$ on $S'$. Given to circle packing metrics $r$ and $R$ on $(S, \mathcal{T}, I)$, we can define corresponding circle packing metrics $r'$ and $R'$ on $(S', \mathcal{T}', I')$. If the premises of the conjecture (ie $r|_B \leq R|_B$ and $K_r|_A \leq K_r|_A$) hold on the disk, they will hold on the doubled surface. Therefore if we find a interior vertex $i \in A$ such that $r(i) > R(i)$, this violation will carry over to the doubled surface and thus refute the conjecture for the closed surface.

**Conjecture 3.1.11** Let $(S, \mathcal{T}, I)$ be a connected weighted triangulated disk in Euclidean geometry such that $I \geq 1$. Let $A$ be the interior vertices and $B$ be the boundary vertices. Then for all circle packing metrics $R$ and $r$ on $(S, \mathcal{T}, I)$ where $R|_B \geq r|_B$ and $K_R|_A \geq K_r|_A$, $R \geq r$.

Thus if `Conjecture 3.10` holds for Euclidean background geometry then `Conjecture 3.1.11` holds for Euclidean background geometry. Additionally if we find a counterexample for `Conjecture 3.1.11`, then we have a counterexample for `Conjecture 3.10` for Euclidean background geometry.

3.2. **Construction of the Counterexample.**

We will now construct a weighted triangulated surface in euclidean geometry with $I \geq 1$ which can be doubled to create a counterexample for the conjecture. Let $V = [4]$, $A = \{4\}$, and $B = [3]$. Thus $A \sqcup B = V$. Let the edges be $E = \{\{1,2\}, \{1,3\}, \{2,3\}, \{1,4\}, \{2,4\}, \{3,4\}\}$ and let the faces be $F = \{\{1,2,4\}, \{1,3,4\}, \{2,3,4\}\}$. Let $I : E \to [1, \infty]$ where

$$I(e) := \begin{cases} 4 \text{ if } e = (2,4) \\ 3 \text{ if } e = (3,4). \\ 1 \text{ otherwise} \end{cases} \tag{53}$$

Thus $(\mathbb{D}, \{V, E, F\}, I)$ is a weighted triangulated disk.

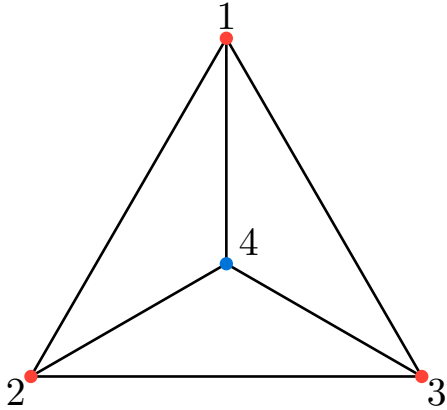

Figure 2. Triangulation of Counter Example

Let $r, R : V \to \mathbb{R}_{>0}$, where

$$r_E = (100, 100, 100, 155) \text{ and } R_E = (110, 240, 220, 150). \tag{54}$$

Note we are using the vector definition of functions from $V$ to $\mathbb{R}$.

### 3.3. Verification of Counterexample.

#### 3.3.1. *Verification of Circle Packing Metrics.*

To verify that such construction forms a valid counterexample, we first need to verify that $r$, $R$ are valid circle packing metrics. To do this we have to calculate the edge lengths using (1). The induced edge lengths for $r$ is

| Edge $e$ | Edge Length $l_r(e)$ |
|---|---|
| $\{1,2\}$ | 200.0 |
| $\{1,3\}$ | 200.0 |
| $\{2,3\}$ | 200.0 |
| $\{1,4\}$ | 255.0 |
| $\{2,4\}$ | $\sqrt{100^2 + 155^2 + 2(155)(100)(4)} \approx 397.52358$ |
| $\{3,4\}$ | $\sqrt{100^2 + 155^2 + 2(155)(100)(3)} \approx 356.40567$ |

and the induced edge lengths for $R$

| Edge $e$ | Edge Length $l_R(e)$ |
|---|---|
| $\{1,2\}$ | 350.0 |
| $\{1,3\}$ | 330.0 |
| $\{2,3\}$ | 460.0 |
| $\{1,4\}$ | 260.0 |
| $\{2,4\}$ | $\sqrt{240^2 + 150^2 + 2(150)(240)(4)} \approx 606.71245$ |
| $\{3,4\}$ | $\sqrt{220^2 + 150^2 + 2(220)(150)(3)} \approx 518.55569$ |

We have to use the edge lengths to verify the triangle inequality. For side lengths $l_1, l_2$, and $l_3$, to see if the lengths $l_1$, $l_2$, and $l_3$ form a triangle it suffices to see if they satisfy the following inequality:

$$(l_1 + l_2 + l_3)(l_1 + l_2 - l_3)(l_1 + l_3 - l_2)(l_2 + l_3 - l_1) > 0 \tag{55}$$

If we plug in the values for $l_r$ and $l_R$ for each triangle we can see they satisfy the triangle inequality. Thus $r$ and $R$ are valid circle packing metrics.

### 3.4. Computation of Discrete Curvature.

To compute the induced discrete curvature at a vertex, we first need to compute the inner angle of each triangle the vertex is a part of. To compute the inner angle at a vertex $a$, we can use the Euclidean cosine law. Let $A$ be the length of edge opposite to $a$ and $B$ and $C$ the length of the other edges.

$$\cos(\theta_a) = \frac{B^2 + C^2 - A^2}{2BC} \Rightarrow \theta_a = \arccos\left(\frac{B^2 + C^2 - A^2}{2BC}\right) \tag{56}$$

Once we calculate the inner angles of each triangle the given vertex is a part of, we can plug the values into the discrete curvature formula (4).

Thus the induced discrete curvatures for $r$ is

| Vertex $e$ | Discrete Curvature $K_r(e)$ (up to 5 decimal places) |
|---|---|
| 1 | 2.37781 |
| 2 | 4.59519 |
| 3 | 4.00207 |
| 4 | 4.73289 |

Thus the induced discrete curvatures for $R$ is

| Vertex $e$ | Discrete Curvature $K_R(e)$ (up to 5 decimal places) |
|---|---|
| 1 | 1.21223 |
| 2 | 5.21346 |
| 3 | 4.51403 |
| 4 | 4.76824 |

## 3.5. Verification of Counterexample.

By the definition of $r$ and $R$, we see that $r|_B = (100, 100, 100) \leq R|_B = (110, 240, 220)$. Additionally, $K_r(4) \leq K_R(4)$, but $r(4) > R(4)$.

## 3.6. Conclusion of the Counter Example.

The premises hold $(r|_B \leq R|_B, K_r|_A \leq K_R|_A)$ but $r(4) > R(4)$. This violates the conjecture on the disk. By a doubling argument, it also provides a counterexample for a closed surface.

## 4. Acknowledgements

The research of this author is supported by the National Science Foundation of the United States of America under Grant #2220271. This work was done with the support of the DIMACS REU at Rutgers University, under the mentorship of Professor Feng Luo, Professor Hongbin Sun, Dr. Zhenghao Rao, and Mr. Kuijin Liu.

## References

1. Rodin, B., Sullivan, D.: The convergence of circle packings to the Riemann mapping. Journal of Differential Geometry. 26, 349–360 (1987). https://doi.org/10.4310/jdg/1214441375
2. Thurston, W.P.: The geometry and topology of three-manifolds: With a preface by Steven P. Kerckhoff. American Mathematical Society (2022)
3. Packing, V., Bowers, P., Stephenson, K.: Uniformizing Dessins And Belyi Maps Via Circle Packing. Mem. Amer. Soc. 170, (1998)
4. Bowers, P.L., Hurdal, M.K.: Planar Conformal Mappings of Piecewise Flat Surfaces. In: Hege, H.-C. and Polthier, K. (eds.) Visualization and Mathematics III. pp. 3–34. Springer Berlin Heidelberg, Berlin, Heidelberg (2003)
5. Stephenson, K., Cannon, J., Floyd, W., Parry, W.: Introduction to circle packing: The theory of discrete analytic functions. The Mathematical Intelligencer. 29, 63–66 (2007). https://doi.org/10.1007/BF02985693
6. Guo, R.: Local rigidity of inversive distance circle packing, https://arxiv.org/abs/0903.1401
7. Stephenson, K.: Circle packing: a mathematical tale. Notices of the AMS. 50, 1376–1388 (2003)

8. Beardon, A.F., Stephenson, K.: The Schwarz-Pick lemma for circle packings. Illinois Journal of Mathematics. 35, 577–606 (1991)
9. Eeuwen, J.V.: The Discrete Schwarz-Pick Lemma for Overlapping Circles. Proceedings of the American Mathematical Society. 121, 1087–1091 (1994)
10. Zhou, Z.: Circle patterns with obtuse exterior intersection angles, https://arxiv.org/abs/1703.01768
11. Verdière, Y.C. de: Un principe variationnel pour les empilements de cercles. Inventiones mathematicae. 104, 655–669. https://doi.org/10.1007/BF01245096
12. Luo, F.: Rigidity of Polyhedral Surfaces, https://arxiv.org/abs/math/0612714
13. Xu, X.: Rigidity of inversive distance circle packings revisited, https://arxiv.org/abs/1705.02714
14. Chow, B., Luo, F.: Combinatorial Ricci Flows on Surfaces, https://arxiv.org/abs/math/0211256

Mathematics Department, University of Texas at Austin, Austin, Texas 75033
*Email address:* arham.lodha@utexas.edu